\documentclass[final,5p,times,twocolumn]{elsarticle}

 \usepackage{graphics}

\usepackage{epsfig}

\usepackage{amsmath}
\usepackage{amsthm}
\usepackage{amsbsy}
\usepackage{esint}
\usepackage{amssymb}

\biboptions{super}

\newcommand{\bsq}{\mathbf{b}^2}
\newcommand{\s}{\mathbf{s}}
\renewcommand{\c}{\mathbf{c}}
\newcommand{\tr}{\textrm{tr}}
\renewcommand{\S}{\mathbf{S}}
\newcommand{\Wi}{\mathrm{Wi}}

\journal{Journal of Non-Newtonian Fluid Mechanics}

\begin{document}

\begin{frontmatter}



\title{Symmetric factorization of the conformation tensor in viscoelastic fluid models}

\author[allofus]{Nusret Balci}
\address[allofus]{Institute for Mathematics and its Applications, University of Minnesota, Minneapolis, MN 55455-0134}

\author[allofus,becca]{Becca Thomases}
\address[becca]{Department of Mathematics, University of California, Davis, CA 95616}

\author[allofus,renardy]{Michael Renardy}
\address[renardy]{Department of Mathematics, Virginia Tech, Blacksburg, VA  24061-0123}

\author[allofus,charlie1,charlie2,charlie3]{Charles R. Doering\corref{*}}
\cortext[*]{E-mail: doering@umich.edu}
\address[charlie1]{Department of Physics, University of Michigan, Ann Arbor, MI 48109-1040}
\address[charlie2]{Department of Mathematics, University of Michigan, Ann Arbor, MI 48109-1043}
\address[charlie3]{Center for the Study of Complex Systems, University of Michigan, Ann Arbor, MI 48109-1107}

\begin{abstract}
The positive definite symmetric polymer conformation tensor possesses a unique symmetric square root that satisfies a closed evolution equation in the Oldroyd-B and FENE-P models of viscoelastic fluid flow.
When expressed in terms of the velocity field and the symmetric square root of the conformation tensor, these models' equations of motion formally constitute an evolution in a Hilbert space with a total energy functional that defines a norm.
Moreover, this formulation is easily implemented in direct numerical simulations resulting in significant practical advantages in terms of both accuracy and stability.
\end{abstract}

\begin{keyword}

\end{keyword}

\end{frontmatter}



\section{Introduction}
Familiar models of viscoelastic polymeric fluids present challenging problems for both mathematical analysis and numerical computations.
One of the difficulties stems from the nature of the stress evolution equations.
Although there is indeed some stress diffusion, the diffusion of polymers is typically orders of magnitude smaller than for non-polymeric molecules and so is often neglected in direct numerical simulations.
The difficulties manifest themselves both in the form of loss of accuracy and stability in numerical schemes, and in the absence of effective {\it a priori} estimates for analysis.
Despite recent progress in the field, many important problems remain open.

In this paper we focus on two models, Oldroyd-B and FENE-P.
It has proven to be a difficult task to devise numerical schemes that are efficient, accurate, and stable at the same time.
One way to ease the numerical problems is to add artificially large stress diffusion, and this has been done for a long time.
Fattal and Kupferman proposed a log-conformation scheme directly evolving the matrix logarithm of the positive definite conformation tensor\cite{Kupferman2005} that, according to their reports, indeed helps with stability issues.
Another method developed by Collins {\it et al} evolves the eigenvalues of the conformation tensor\cite{Collins2006}.
Lozinski and Owen also proposed to work with the deformation tensor\cite{Lozinski2003}, another of the square roots of the conformation tensor.

We consider a square root method as well, but unlike any previous work we are aware of, we derive an evolution equation for the {\em positive-definite} square root by taking advantage of the $O(n)$ degeneracy in the matrix square root in $n$ dimensions.
This turns out to be both theoretically and numerically convenient.
On the one hand it allows the dependent variables in the Oldroyd-B and FENE-P models to take values in a vector space with a natural norm defined by the physical energy.
On the other hand we observe that, at practically no additional computational cost, this formulation produces significant gains in both numerical stability and numerical accuracy---without adding any artificial stress diffusion---as compared to directly evolving the conformation tensor. 
We note in particular that the flow studied in this paper has a hyperbolic stagnation point and our simulations for the Oldroyd-B model extend far beyond the Weissenberg number at which the stress in the associated steady flow becomes infinite.

\section{Mathematical Framework}

The nondimensional equations of motion are
\begin{equation}
\frac{\partial\mathbf{u}(\mathbf{x},t)}{\partial t}+\mathbf{u}\cdot\nabla\mathbf{u}+\nabla p = \frac{1}{\mathrm{Re}}\Delta\mathbf{u}+\nabla\cdot\boldsymbol{\tau} + \mathbf{f}(\mathbf{x}, t),\quad\nabla\cdot\mathbf{u}=0,\label{NS}
\end{equation}
with $\mathbf{x} \in \mathcal{R}^n$ ($n=2$ or $3$) and Reynolds number $\mathrm{Re}=U\ell/\nu$ where $U$ and $\ell$ represent appropriate choices of velocity and length scales for the problem under investigation.
The externally applied body force is denoted by $\mathbf{f}$, and the polymer stress tensor $\boldsymbol{\tau}(\mathbf{x},t)$ is
\begin{equation}
\boldsymbol{\tau}= - \frac{s}{\mathrm{Re}}\, \mathbf{s}(\mathbf{c})
\end{equation}
where the symmetric positive definite polymer conformation (a.k.a. configuration) tensor $\mathbf{c}(\mathbf{x},t)$ evolves according to
\begin{equation}
\frac{\partial\mathbf{c}}{\partial t}+\mathbf{u}\cdot\nabla\mathbf{c}=\mathbf{c}\nabla\mathbf{u}+(\nabla\mathbf{u})^{T}\mathbf{c}+\mathbf{s}(\mathbf{c}). \label{eq:configeqn}
\end{equation}
The parameter $s$ is a coupling constant proportional to the concentration of the polymers in the fluid, and the tensor $\mathbf{s}(\mathbf{c})$ takes different forms in various non-Newtonian models.
For the Oldroyd-B model
\begin{equation}\label{OB}
\mathbf{s}(\mathbf{c})=\frac{1}{\mathrm{Wi}}(\mathbf{I}-\mathbf{c})
\end{equation}
where the Weissenberg number $\mathrm{Wi}=U\lambda/\ell$ is the product of the polymer relaxation time $\lambda$ and the rate of strain $U/\ell$.
For the FENE-P model
\begin{equation}\label{FP}
\mathbf{s}(\mathbf{c})=\frac{1}{\mathrm{Wi}}\left( \mathbf{I}-\frac{\mathbf{c}}{1-(\mathrm{tr}\mathbf{c}/l^2)}\right),
\end{equation}
where $l^2$ is proportional to the square of the maximum polymer length.
For both models the total mechanical energy of the system is the sum of the fluid's kinetic energy and the elastic potential energy of the polymers:
\begin{equation}
\mathcal{E}(t)=\frac{1}{2}\int\left[|\mathbf{u}(\mathbf{x},t)|^{2}+\mathrm{tr} \, \boldsymbol{\tau}\right]dx\, dy\, dz. \label{energy}
\end{equation}
This energy\footnote{$\mathcal{E}$ does not include the entropic term that contributes to the free energy of the system.} is formally conserved by the dynamics in the limits of infinite Reynolds and Weissenberg numbers.

Unlike the situation for Newtonian fluids modeled by the incompressible Navier-Stokes (or Stokes) equations,
the total energy does not define a natural norm, or even a metric, in the phase space of the dynamical fields $\mathbf{u}$ and $\mathbf{c}$.
Indeed, the $(\mathbf{u},\mathbf{c})$ phase space is not even a linear vector space.
This mathematical awkwardness results from the fact that the relevant space for the conformation tensor $\mathbf{c}$, the space of symmetric positive definite matrices, is not a vector space: linear combinations of positive matrices are not necessarily positive.
These facts complicate the analysis of these models and preclude implementation of useful techniques including nonlinear (energy) stability notions\cite{DES2006}.

This problem can be circumvented, however, by reformulating the models in terms of the (unique) symmetric square root  $\mathbf{b}(\mathbf{x},t)$ of the conformation tensor $\mathbf{c}(\mathbf{x},t)$.
We write
\begin{equation}
c_{ij}(\mathbf{x},t) = \sum_{k=1}^n b_{ik}(\mathbf{x},t) b_{kj}(\mathbf{x},t) \quad \text{with} \quad b_{ij}(\mathbf{x},t)=b_{ji}(\mathbf{x},t),
\end{equation}
so the polymer energy density is a function of the matrix norm of $\mathbf{b}$,
\begin{equation}
\|\mathbf{b}\|^2 = \sum_{i,j=1}^{n} b_{ij}^{2} = \mathrm{tr}(\mathbf{b}^{T}\mathbf{b}) = \mathrm{tr} \, \mathbf{c}.
\end{equation}
The work required to implement this proposal is to precisely articulate the dynamics of $\mathbf{b}$, a not altogether trivial task due to the inherent degeneracy of the matrix square root.

In the Oldroyd-B case solutions of differential equations of the form
\begin{equation}
\left(\frac{\partial}{\partial t}+\mathbf{u}\cdot\nabla\right)\mathbf{b}=\mathbf{b}\nabla\mathbf{u}+\mathbf{a}\mathbf{b}+\frac{1}{2\mathrm{Wi}}((\mathbf{b}^{T})^{-1}-\mathbf{b}),\label{eq:Bevol}
\end{equation}
where $\mathbf{a}(\mathbf{x},t)$ is {\em any} antisymmetric matrix, satisfy $\mathbf{b}^{T}\mathbf{b}=\mathbf{c}$ pointwise in space and time when the initial data satisfy $\mathbf{b}^{T}(\mathbf{x},0)\mathbf{b}(\mathbf{x},0)=\mathbf{c}(\mathbf{x},0)$.
Likewise, in the FENE-P case the evolution
\begin{equation}
\left(\frac{\partial}{\partial t}+\mathbf{u}\cdot\nabla\right)\mathbf{b}=\mathbf{b}\nabla\mathbf{u}+\mathbf{a}\mathbf{b}+\frac{1}{2\mathrm{Wi}}\left((\mathbf{b}^{T})^{-1}-\frac{\mathbf{b}}{1-\|\mathbf{b}\|^{2}/l^2}\right),\label{eq:BevolFENEP}
\end{equation}
produces such a square root of $\mathbf{c}$ when $\mathbf{a}(\mathbf{x},t)$ is any antisymmetric matrix and $\mathbf{b}^{T}(\mathbf{x},0)\mathbf{b}(\mathbf{x},0)=\mathbf{c}(\mathbf{x},0)$.
The key observation is that by choosing $\mathbf{a}(\mathbf{x},t)$ properly we can tune the evolution equations (\ref{eq:Bevol}) and (\ref{eq:BevolFENEP})---and similar models with an upper convective derivative---to preserve the symmetry of $\mathbf{b}$.
That is, starting with symmetric initial data $\mathbf{b}^{T}(\mathbf{x},0)=\mathbf{b}(\mathbf{x},0)$, the subsequent evolution will preserve the symmetry.

Toward this end we write $\nabla\mathbf{u}=(u_{ij})=(u_{j,i})$ and
\begin{equation}\label{matrix:a}
\mathbf{a}=\left(\begin{array}{ccc}
0 & a_{12} & a_{13}\\
-a_{12} & 0 & a_{23}\\
-a_{13} & -a_{23} & 0
\end{array}\right),\quad i,j=1,2,3,
\end{equation}
in $n=3$ spatial dimensions and
\begin{equation}
\mathbf{a}=\left(\begin{array}{cc}
0 & a_{12}\\
-a_{12} & 0
\end{array}\right)
\end{equation}
in  $n=2$ dimensions.
Define
\begin{equation}
\mathbf{r}=(r_{ij})=\mathbf{b} (\nabla\mathbf{u})+\mathbf{a}\mathbf{b}.
\end{equation}
We now show that we may choose the matrices $\mathbf{a}$, depending on $\nabla\mathbf{u}$ and the symmetric matrix $\mathbf{b}$ pointwise in space and time, so that $\mathbf{r}$ is a field of symmetric matrices, i.e., $r_{ij}=r_{ji}$.
For $n=3$ the explicit formulas for the elements $a_{ij}$ come from solving the system of $3$ linear equations
\begin{eqnarray}
(b_{11}+b_{22})a_{12}+b_{23}a_{13}-b_{31}a_{23} &=& w_{1},\\
b_{23}a_{12}+(b_{11}+b_{33})a_{13}+b_{12}a_{23} &=& w_{2},\\
-b_{13}a_{12}+b_{12}a_{13}+(b_{22}+b_{33})a_{23} &=& w_{3},
\end{eqnarray}
where
\begin{eqnarray}
w_{1} \ = \ (b_{12}u_{1,1}-b_{11}u_{2,1}) &+& (b_{22}u_{1,2}-b_{21}u_{2,2}) \nonumber \\
&+& (b_{32}u_{1,3}-b_{31}u_{2,3}),\label{w:1}\\
w_{2} \ = \ (b_{13}u_{1,1}-b_{11}u_{3,1}) &+& (b_{33}u_{1,3}-b_{31}u_{3,3}) \nonumber \\
&+& (b_{23}u_{1,2}-b_{21}u_{3,2}),\label{w:2}\\
w_{3} \ = \ (b_{13}u_{2,1}-b_{12}u_{3,1}) &+& (b_{23}u_{2,2}-b_{22}u_{3,2}) \nonumber \\
&+& (b_{33}u_{2,3}-b_{32}u_{3,3})\label{w:3}.
\end{eqnarray}
In matrix notation, this is the system of equations
\begin{equation}
\left(\begin{array}{ccc}
b_{11}+b_{22} & b_{23} & -b_{31}\\
b_{23} & b_{11}+b_{33} & b_{12}\\
-b_{31} & b_{12} & b_{22}+b_{33}
\end{array}\right)\left(\begin{array}{c}
a_{12}\\
a_{13}\\
a_{23}
\end{array}\right)=\left(\begin{array}{c}
w_{1}\\
w_{2}\\
w_{3}
\end{array}\right).
\end{equation}
Then by swapping the first and the third columns of the coefficient matrix (and hence, also $a_{23}$ and $a_{12}$), and subsequently swapping the first and the third rows of the resulting coefficient matrix (and hence, also $w_{1}$ and $w_{3}$), and finally multiplying the second row and the second column of the resulting matrix by $-1$ (and hence, also replacing $a_{13}$ and $w_{2}$ by $-a_{13}$ and $-w_{2}$, respectively), we obtain
\begin{equation}
\left( \mathbf{tr}(\mathbf{b})\mathbf{I}-\mathbf{b} \right) \, \tilde{\mathbf{a}} = \mathbf{v},
\end{equation}
where
\begin{equation}
\tilde{\mathbf{a}}=\left(\begin{array}{c}
a_{23}\\
-a_{13}\\
a_{12}
\end{array}\right),\quad\mathbf{v}=\left(\begin{array}{c}
w_{3}\\
-w_{2}\\
w_{1}
\end{array}\right).
\end{equation}
When $\mathbf{b}$ is symmetric at the point $(\mathbf{x},t)$ there is an orthogonal matrix $\mathbf{p}(\mathbf{x},t)$ such that
\begin{equation}
\mathbf{b}=\mathbf{p}^{T}\mathrm{diag}\left\{ \lambda_{1},\lambda_{2},\lambda_{3}\right\} \mathbf{p},
\end{equation}
where $\lambda_{i}$ are the eigenvalues of $\mathbf{b}$.
Thus, we have
\begin{align}
\mathrm{trace}(\mathbf{b})\mathbf{I}-\mathbf{b} & =\mathrm{trace}(\mathbf{b})\mathbf{I}-\mathbf{p}^{T}\mathrm{diag}\left\{ \lambda_{1},\lambda_{2},\lambda_{3}\right\} \mathbf{p}\\
 & =\mathbf{p}^{T}(\mathrm{trace}(\mathbf{b})\mathbf{I}-\mathrm{diag}\left\{ \lambda_{1},\lambda_{2},\lambda_{3}\right\} )\mathbf{p}\\
 & =\mathbf{p}^{T}\mathrm{diag}\left\{ \lambda_{2}+\lambda_{3},\lambda_{1}+\lambda_{3},\lambda_{1}+\lambda_{2}\right\} \mathbf{p}
\end{align}
Again, assuming that $\mathbf{b}$ is positive definite (although this condition can be clearly relaxed to include a large class of semidefinite objects) we can solve for $\mathbf{a}$ uniquely so that the evolution equations (\ref{eq:Bevol}) and (\ref{eq:BevolFENEP}) used to obtain $\mathbf{b}$ at later times are symmetrized.
The explicit algebraic formulas for the elements $a_{ij}$ for $n=3$ are displayed in the appendix.
In the much simpler case of $n=2$ space dimensions we have
\begin{equation}
a_{12}=\frac{(b_{12}u_{1,1}-b_{11}u_{2,1})+(b_{22}u_{1,2}-b_{21}u_{2,2})}{b_{11}+b_{22}}.
\end{equation}

This construction puts the full dynamics back in a vector space setting (the direct product of vector fields $\mathbf{u}$ and symmetric tensor fields $\mathbf{b}$).
For Oldroyd-B the energy functional (\ref{energy}) is proportional to the vector norm (squared) on the direct product space (modulo an additive constant).
For the FENE-P model, the convexity of $ \mathbf{b}^{T}\mathbf{b}/(1-(\text{tr}(\mathbf{b}^{T}\mathbf{b})/l^2))$ in a neighborhood of the origin (i.e., for $\|\mathbf{b}\| < l$) allows as well for a natural energy norm.
And as shown in the next section, in some cases this reformulation of the dynamics leads to significant practical improvements in direct numerical simulations.

\section{Numerical Experiments}

\begin{figure}[htp]
\centering
\includegraphics[width=3.5in]{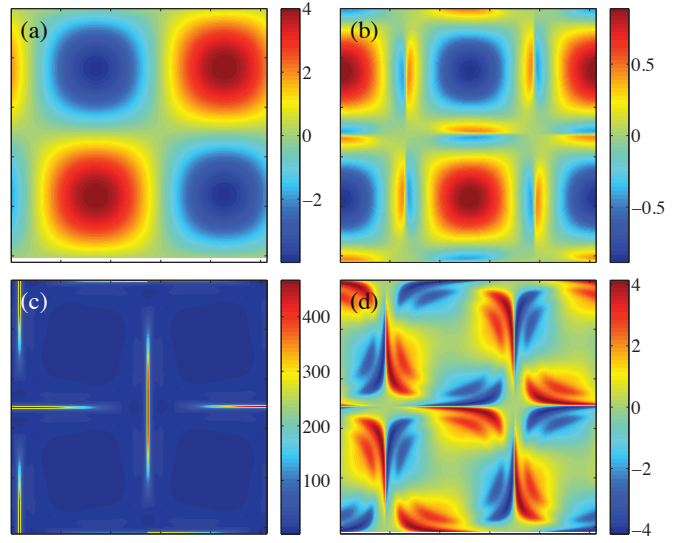}
\caption{(a) Contour plots of $\textrm{curl}\mathbf{f}$ for the force given by Eq. \eqref{force}. (b) - (d) Contour plot of vorticity, $\tr\mathbf{c},$ and $\mathbf{c}_{12}$ 
for isotropic initial data, $\Wi = 5$ at $t=10$.   }
\label{forcefig}
\end{figure}

As a test of the numerical accuracy and stability of the square-root method we consider the zero-Reynolds number (Stokes) limit of the Oldroyd-B and FENE-P models in which case the momentum equation \eqref{NS} reduces to
\begin{align}
\nabla p & =\Delta\mathbf{u}+\nabla\cdot\boldsymbol{\tau}+\mathbf{f},\quad\nabla\cdot\mathbf{u}=0.\label{Stokes}\end{align}
Here $\tau=-s\s(\c)$ with $\s(\c)$ given by \eqref{OB} or \eqref{FP} for the Oldroyd-B model or FENE-P model respectively.
In what follows we fix $s=0.5.$ Following recent studies\cite{TS2007, TS2009} we consider a $2\pi-$periodic domain in $n=2$ space dimensions ($[-\pi,\pi]^2$) and impose a steady background force
\begin{equation}
\label{force}\mathbf{f}=(-2\sin x\cos y, 2\cos x \sin y),
\end{equation}  $\textrm{curl}\mathbf{f}$ is shown in Fig. \ref{forcefig}(a).  
In the absence of polymer stress this yields a four-roll mill geometry for the fluid velocity.
One well-known consequence of this body-force imposed extensional geometry in the Oldroyd-B model is that the polymer
stress and stress gradients grow exponentially in time \cite{TS2007, Renardy2006, RH1988} and inevitably produce numerical problems.
In particular, when resolving steep gradients the loss of positive definiteness of the conformation tensor due to numerical error can lead to breakdown of the computational schemes.
One common solution to these difficulties is the addition of artificial polymer stress diffusion.
Although some polymer stress diffusion can be derived from the basic physics in the model, the magnitude of the physically relevant diffusion is far too small to have an effect on numerical simulations\cite{KL1989}.
In the following we do not add {\it any} stress diffusion to the numerical calculations.

The Stokes-Oldroyd-B system (and FENE-P) is solved with a pseudo-spectral method\cite{Peyret2002}.
In the usual formulation the conformation tensor is evolved using a second-order Adams-Bashforth method.
The initial data is prescribed, and given $\c$ the Stokes equation is inverted in Fourier space for $\mathbf{u}$.
Given $\mathbf{u}$ the nonlinearities of the conformation tensor evolution are evaluated using a smooth filter applied in Fourier space before the quadratic terms are multiplied in real space\cite{HL2007}.
The conformation tensor is then discretized on the Fourier transform side.  It should be noted that the numerical implementation of the square-root method adds no computational cost.  

In a recent investigation\cite{TS2007} the Stokes-Oldroyd-B equations were solved starting from isotropic, i.e., $\c(0)=\mathbf{I},$ initial data and the stress was observed to diverge (exponentially in time) at the extensional stagnation points in the flow for sufficiently large Weissenberg number.
However, outside of an exponentially decreasing region around the extensional stagnation point, the solutions became steady after an initial transient.
These near-steady solutions preserve many symmetries: the stress is symmetric and aligned along the direction of extension and the flow has an underlying four-roll structure.
For sufficiently large $\Wi$ additional oppositely signed vortices arise along the stable and unstable manifolds of the extensional stagnation point. Figure
\ref{forcefig}(b) - (d) displays these symmetric solutions in the case $\Wi = 5,$ at $t=10.$  
These symmetries are broken as the initial data is perturbed and it was shown that instabilities arise for sufficiently large $\Wi$\cite{TS2009,TST2010}.

In what follows we discuss both accuracy and stability improvements for the Oldroyd-B and FENE-P models using the square root method.  In subsection \ref{sec:accuracy} we consider homogeneous isotropic initial data $\c(\mathbf{x},0)=\bsq(\mathbf{x},0)=\mathbf{I}$ and compare solutions to the Stokes-Oldroyd-B model obtained using the direct evolution of $\c$ with those obtained by evolving the symmetric square root.
In subsection \ref{sec:stability} we consider perturbations to the initial data for the Stokes-Oldroyd-B model (as in previous studies\cite{TS2009, TST2010}) to see how far the simulations run, for a fixed resolution using each method, before numerical divergences appear.
Finally, in subsection \ref{sec:fenep} we revisit both of these questions for the FENE-P model.

\subsection{Accuracy}\label{sec:accuracy}

\begin{figure}[htp]
\centering
\includegraphics[width=3.5in]{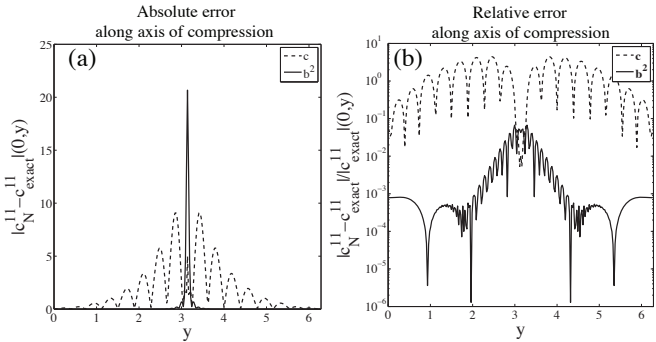}
\caption{(a) Absolute error $|\c_N-\c_{\textrm{exact}}|$ measured along the axis of compression of the first component of the
conformation tensor with $N^2=256^2,$ $\Wi = 5$, at $t = 10.$  Dotted line is $\c$ and
solid line is $\bsq.$  (b) Relative error $\frac{|\c_N-\c_{\textrm{exact}}|}{|c_{\textrm{exact}}|}.$  }
\label{Relerr}
\end{figure}

Figure \ref{Relerr} shows the difference between the solution to Stokes-Oldroyd-B with $N^2=256^2$ grid cells and the ``exact'' solution computed by evolving $\c$ with $N^2=2048^2$ grid cells, resolved to at least 6 digits of accuracy.
The dotted line shows the solution computed by directly evolving the conformation tensor $\c$, while the solid line shows the solution computed by evolving the symmetric square root.
The simulation is performed with $\Wi=5$ and the result of the computation is shown at $t=10.$
Panel (a) shows the absolute error ($|\c_N-\c_{\textrm{exact}}|$) in the first component of the conformation tensor ($\c_{11}$) along the direction of compression because this is precisely where steep gradients form.
We observe that away from the extensional stagnation point the square root gives a much better approximation and it is only very near the extensional stagnation point where the evolution of $\c$ gives a better approximation.
The relative error is shown in panel (b) to emphasize that although the square-root method's error at the extensional stagnation point is larger, the relative error is actually quite small because the stress is very large there.

\begin{figure}[htp]
\centering
\includegraphics[width=3.5in]{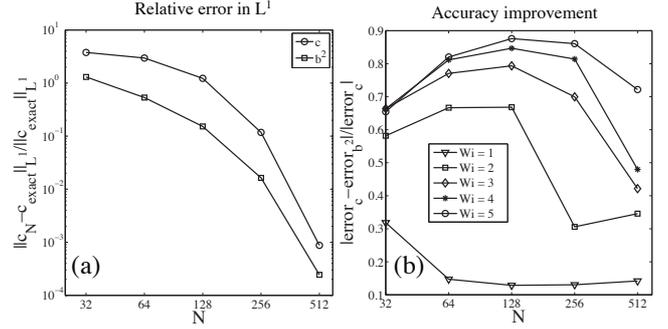}
\caption{(a) Relative error in $L^1$ shows $\displaystyle \frac{\|\c_{N}-\mathbf{c_{\textrm{exact}}}\|_{L^1}}{\|\c_{\textrm{exact}}\|_{L^1}}$ computed with $\c$ and $\bsq$, for $\Wi = 5$ at $t=10.$
(b) Improvement in accuracy as a function of $N,$ for $\Wi =1 - 5$ at $T=t/\Wi=2.$ }
\label{accuracy}
\end{figure}

Figure \ref{accuracy} (a) shows the relative error
\begin{equation}
\frac{\|\c_{N}-\mathbf{c_{\textrm{exact}}}\|_{L^1}}{\|\c_{\textrm{exact}}\|_{L^1}}
\end{equation}
measured in the $L^1([-\pi,\pi]^2)$ norm for both $\c_N$ and $(\bsq)_N$ for $N^2=32^2,64^2,128^2,256^2,512^2,$ comparing to the ``exact" solution defined by the $N^2=2048^2$ simulation.
The $L^1$ norm is chosen because it takes into account the average error over the entire domain.
This computation is also for $\Wi = 5$ at $t = 10.$

In this averaged sense we see that the error is always smaller using the square root method.
The improvement in accuracy (in the $L^1-$sense) using the square root method is shown in figure \ref{accuracy} (b).
Here we plot
\begin{equation}
\frac{|\textrm{error}_{\c}-\textrm{error}_{\bsq}|}{|\textrm{error}_{\c}|}
\end{equation}
for $\Wi = 1,2,3,4,5$ for $N^2=32^2,64^2,128^2,256^2,512^2.$  Each simulation is computed at $T = t/\Wi = 2,$ and this scaled-time is used because the solutions grow exponentially like $e^{t/\Wi}.$
There is a significant improvement for higher $\Wi$, in particular for lower resolutions $N^2=128^2$ and $256^2.$

\subsection{Stability}\label{sec:stability}

\begin{figure}[htp]
\centering
\includegraphics[width=3.5in]{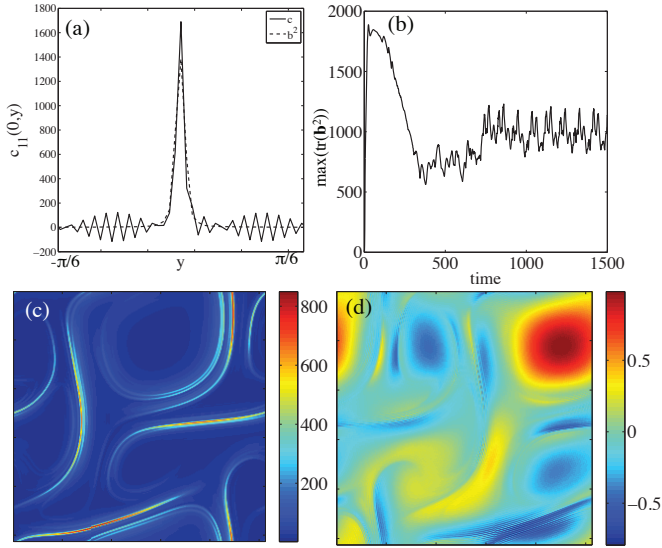}
\caption{(a) Plot of first component of conformation tensor $\c_{11}(0,y)$ and $\bsq_{11}(0,y)$ along direction of compression near the extensional stagnation point
for $\Wi = 10$ at $t=15.$  Computations of $\c$ stop producing finite values at $t=19.$  (b) Plot of $\textrm{max}(\textrm{tr}(\bsq))$ as a function of
time over $0<t<1500.$  (c) $\Wi = 10$, contour plot of $\textrm{tr}(\bsq)$ at $t=1000$ on $[-\pi,\pi]^2.$
(d) $\Wi = 10$, contour plot of the vorticity of the flow field at $t=1000$ on $[-\pi,-\pi]^2.$  All simulations performed with $N^2=256^2$ grid points. }
\label{stability}
\end{figure}

Experiments on low-Reynolds number viscoelastic turbulence \cite{GS2000,GS2001,GS2004} and instabilities in extensional geometries \cite{ATDG2006} have inspired many numerical studies of low-Reynolds number viscoelastic fluids\cite{PAO2007, BBBCM2008,XG2008, TS2009, TST2010}.
Two main instabilities were observed in one investigation\cite{TS2009, TST2010}:  first, for sufficiently large $\Wi$, if a small perturbation is introduced in the initial stress conformation the extensional stagnation point in the flow becomes unstable and loses the pinning to the background steady force.
For larger $\Wi$ other stagnation points lose their pinning to the background force and higher oscillations arise in the flow.
These instabilities occur on long time scales and some artificial polymer stress diffusion was introduced to fully resolve the stress \cite{TS2009, TST2010}.
Here we test these same perturbations to the initial data and for a fixed resolution $N^2=256^2$ and run the same simulations {\it without} any artificial polymer stress diffusion.

Figure \ref{stability} (a) shows a plot of the first component of the conformation tensor for $\Wi = 10$ at $t=15$ computed both by evolving the conformation tensor (solid line) and the square root (dotted line).
The  plot is shown along the axis of compression and it is evident that the stress has accumulated significantly and is quite large near the extensional stagnation point in the flow ($y=0$).
The oscillations produced in the direct evolution of $\c$ lead to the loss of positive-definiteness of the conformation tensor, and the numerical scheme breaks down.
This figure shown is at $t=15$ and the computation fails to produce finite numbers at $t=20$.
However using the square root one can run these simulations to $t=1500$ and even beyond.
Figure \ref{stability} (b) shows a plot of $\max(\tr(\bsq))$ as a function of time for $0<t<1500$.
It is important to point out that although $\max(\tr(\bsq))$ remains bounded in this case (with {\it no} artifical stress diffusion) this quantity clearly depends on $N$ and this is one way that the accuracy of the solution is lost.
However, it is not clear that this level of loss of accuracy is entirely relevant to the flow because the region where $\tr(\bsq)$ gets large diminishes exponentially in time even as $\tr(\bsq)$ grows exponentially in time\cite{TS2007}.
Figure \ref{stability} (c) and (d) show contour plots of $\tr(\bsq)$ and the vorticity of the flow on $[-\pi,\pi]^2$.
Here we see that the previously observed instabilities\cite{TS2009, TST2010} are at least qualitatively reproduced.
The time-dependent behavior is similar, too: the four-roll mill structure of the background force is preserved initially, the extensional stagnation point leaves the origin, and eventually time-dependent oscillations arise in the flow.

Of course with fixed resolution and no stress diffusion there is an inevitable loss of accuracy.
This can be seen in Fig. \ref{stability} (c) and (d) in the slightly fuzzy images indicating oscillations while attempting to resolve the steep gradients in the conformation tensor and vorticity.
It is noteworthy that these simulations are performed with no artificial stress diffusion but nevertheless qualitatively reproduce the well-resolved results that utilized artificial diffusion\cite{TS2009,TST2010}.
The same computations simply cannot be performed with a direct evolution of the conformation tensor (in this particular code).
The square root method allows simulations  to run much longer and at much higher Weissenberg number than evolving $\c$ directly allows.
This indicates that it might be possible to use much smaller---closer to the physically realistic quantity---stress diffusion and still obtain reasonably accurate results, although this will not be pursued in this paper.

\subsection{FENE-P}\label{sec:fenep}

\begin{figure}[htp]
\centering
\includegraphics[width=3.5in]{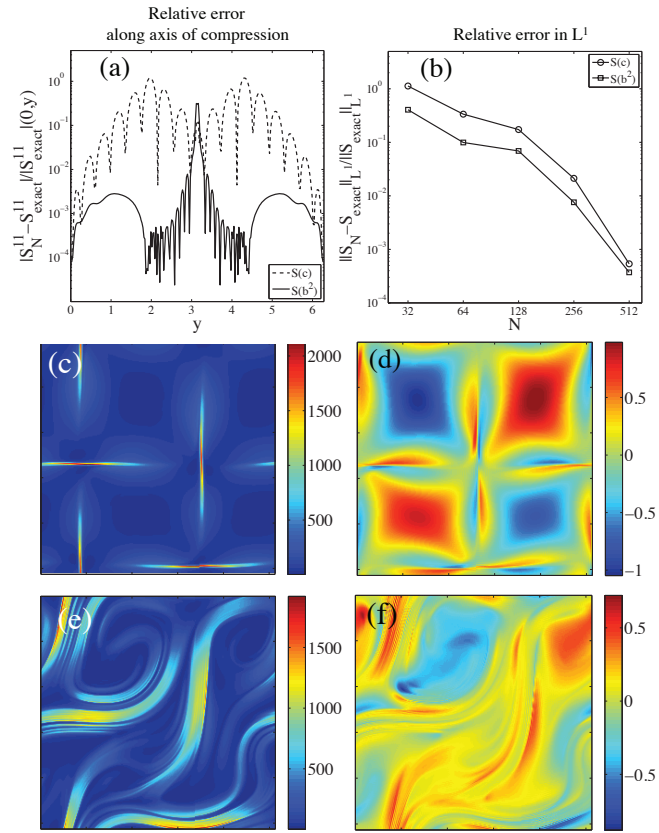}
\caption{(a) Relative error $\frac{|\S_N-\S_{\textrm{exact}}|}{|\S_{\textrm{exact}}|}$ measured along the axis of compression of the first component of $\S(\c)$ for FENE-P with $N^2=256^2,$ $\Wi = 5$ at $t = 10,$ $l^2=100.$  Dotted line is $\S(\c)$ and
solid line is $\S(\bsq).$  (b) Relative error in $L^1$ shows $\displaystyle \frac{\|\S_{N}-\S_{\textrm{exact}}\|_{L^1}}{\|\S_{\textrm{exact}}\|_{L^1}}$ computed with $\S(\c)$ and $\S(\bsq)$, for $\Wi = 5$ at $t=10,$ $l^2=100.$
 (c) $\Wi = 20$, contour plot of $\textrm{tr}(\S(\bsq))$ at $t=100$ on $[-\pi,\pi]^2,$ $l^2=225.$
(d) $\Wi = 20$, contour plot of the vorticity of the flow field at $t=100$ on $[-\pi,-\pi]^2,$ $l^2=225.$  All simulations done with $N^2=256^2$ grid points. (e) $\Wi = 50$, contour plot of $\textrm{tr}(\S(\bsq))$ at $t=500$ on $[-\pi,\pi]^2,$ $l^2=225.$
(f) $\Wi = 50$, contour plot of the vorticity of the flow field at $t=500$ on $[-\pi,-\pi]^2,$ $l^2=225.$  All simulations done with $N^2=256^2$ grid points. }
\label{fenep}\end{figure}

The FENE-P model enforces a limit ($l^2$) on the magnitude of $\tr \, \c$ so the conformation tensor remains bounded.
Steep gradients can still arise in the polymer stress, however, and numerical difficulties remain.
Therefore we also simulated the FENE-P model in a Stokesian solvent to check for possible accuracy and stability improvements by evolving the symmetric square root.

Figure \ref{fenep} (a) and (b) are analogs of Figures \ref{Relerr} (b) and \ref{accuracy} (a) for FENE-P.
The simulations were performed with $\Wi = 5$ and cut-off $l^2=100,$ and are displayed at $t=10$.
Rather than plot the conformation tensor $\c$ and $\bsq$, however, it is more analogous to plot $\S=\frac{\c}{1-(\tr\c/l^2)}$ because this is closely related to the physical stress tensor and includes the factor that gets very large as $\tr\c$ gets near the cut-off $l^2$.
The accuracy improvement is not as large here as it was for Oldroyd-B: for $\Wi=5$ the improvement is about $65\%$ at $N=256^2$ and is only $31\%$ for $N^2=512^2$ but there is still some improvement (especially away from the extensional stagnation point).
As before, the ``exact" solution here comes from a simulation with $N^2=1024^2.$

The significance of the symmetric square root method for FENE-P is much more apparent in terms of stability.
The fact is that we can increase $\Wi$ much more than we can for Stokes-Oldroyd-B. We show results from two different simulations to demonstrate
this.  First in Fig. \ref{fenep} (c) and (d) we show results from perturbed initial data with $\Wi=20$ at $t=100,$ with $l^2=225$. This is just after
the onset of the instability and the flow is still nearly symmetric.  The stress has accumulated along the incoming and outgoing streamlines of the extensional
stagnation point and the four-roll mill structure of the vorticity is still largely preserved. Fig. \ref{fenep} (e) and (f) show results from perturbed initial data with $\Wi=50$ at $t=500,$ with $l^2=225$. 
These same computations evolving $\c$ fail to produce finite values before $t=60$ whereas the evolution of $\mathbf{b}$ appears to continue indefinitely---although, again, there must be some loss of accuracy.
The qualitative behavior is similar to the solutions of Stokes-Oldroyd-B and the instabilities discussed for that case also occur here.
In Fig. \ref{fenep} (c) and (e) we show contour plots of $\tr(\S(\bsq))$ after the instabilities have developed and observe that $\max\tr(\S)$ is quite large.
The time-dependent behavior is also quite complicated and as one can see in Fig. \ref{fenep} (f), the vorticity of the flow is also very complex with many additional vortices continually arising and being destroyed in the flow.

\section{Discussion and Conclusions}

In hindsight both the symmetrization procedure and directly computing the square root evolution equations (\ref{eq:Bevol}) and (\ref{eq:BevolFENEP}) might have been expected to contribute to the gains in stability and accuracy.
Taking the square root reduces large amplitudes which, not unexpectedly, reduces the stiffness in time stepping.\footnote{Clearly, the positive $(2k)^{th}$-roots of $\c$ have entries with even smaller amplitude compared to $\c$ when stress gets very large, which may help the numerics.
This aspect is most vigorously pursued in methods where one computes the logarithm of the matrix $\c$ but these methods can be computationally expensive and more complicated to implement\cite{Kupferman2005}. A comparison of the square-root method with the logarithm method\cite{Kupferman2005} and the method of evolving eigenvalues\cite{Collins2006} is planned for a future study. }  Moreover symmetrizing the system may reduce the stiffness of the time marching as compared to taking $\mathbf{a}=\mathbf{0}$ and simply computing the deformation tensor because components of the symmetric square root matrix will generally have less variance than a square root with no symmetry.
And computing the square root instead of $\c$ has other advantages: the square root computation ensures positivity of $\c$ in the numerical scheme compared to the most direct evolution of the conformation tensor.
We have observed that in practice the square root method can be applied at higher $\Wi$ and for longer time without any artificial numerical stress diffusion than evolving the conformation tensor directly can, enabling one to obtain numerical solutions in more situations.  

One less obvious but perhaps important advantage is the following.
Assuming the locality of modal interactions, i.e., that lower spectral modes of $\c$ are determined predominantly by lower modes of the square root matrix $\mathbf{b}$, we can expect good information about up to the first $2N$ modes of the conformation tensor $\c$ when we know just the first $N$ modes of the square root $\mathbf{b}$.
This speculation basically boils down to the assumption that the Galerkin approximation method works well for both $\mathbf{b}$ and $\c$ for sufficiently large $N$.
Thus we might expect that evolution of the square root improves both stability and accuracy, at least in spectral or pseudospectral schemes.
Whether this is the case in other types of spatial discretizations requires further investigation.
Of course it also remains to implement the full $n=3$ dimensional symmetric square root algorithm and systematically compare its performance with conventional schemes used to investigate, for example, turbulent drag reduction\cite{Graham2010}.

We emphasize that all the numerical simulations shown have been performed without artificial diffusion.  One can add artificial diffusion to the advection of $\mathbf{b},$ however for this to match the standard diffusion of the conformation tensor, a more complicated nonlinear diffusion term is needed.  The square-root method still has limitations for sufficiently large $\Wi$ so it is natural to ask if one could use a more physically realistic diffusion coefficient with the square-root method and this will be pursued in future work. 

The advantage of expressing the Oldroyd-B and FENE-P models as evolutions in a vector space (indeed, a Hilbert space) remains one of theoretical elegance at this point.
Whether or not this formulation might assist the rigorous mathematical analysis of these models is an open question.

\section{Acknowledgements}
The authors are grateful to Michael Graham and Bruno Eckhardt for
enlightening and encouraging discussions.  This work was supported
in part by NSF Awards DMS-0757813 (BT), DMS-0707727 (MR), and
PHY-0855335 (CRD).

\appendix
\section{Entries of the antisymmetric matrix for $n=3$}
The entries $a_{12},a_{13},a_{23}$ of the antisymmetric matrix $\mathbf{a}$ in (\ref{matrix:a}) are given by
\begin{eqnarray}
D \, a_{12} &=& \left(T_{1}T_{2}-B_{3}^{2}\right)w_{1} \ - \ \left(B_{1}T_{1}+B_{3}B_{2}\right)w_{2} \nonumber \\
&\quad& \quad \quad + \  \left(B_{2}T_{2}+B_{1}B_{3}\right)w_{3},\\
D \, a_{13} &=& -\left(B_{1}T_{1}+B_{3}B_{2}\right)w_{1} \ + \ \left(T_{1}T_{3}-B_{2}^{2}\right)w_{2} \nonumber \\
&\quad& \quad \quad - \  \left(B_{2}B_{1}+B_{3}T_{3}\right)w_{3},\\
D \, a_{23} &=& \left(B_{2}T_{2}+B_{1}B_{3}\right)w_{1} \ - \ \left(B_{2}B_{1}+B_{3}T_{3}\right)w_{2} \nonumber \\
&\quad& \quad \quad + \  \left(T_{2}T_{3}-B_{1}^{2}\right)w_{3},
\end{eqnarray}
where
\begin{align}
&D \ \equiv \ \det((\mathrm{tr}\,\mathbf{b})\mathbf{I}-\mathbf{b}) \nonumber \\
&\quad = \ T_{1}\left(T_{2}T_{3} - B_{1}^{2}\right) - B_{2}\left(B_{2}T_{2}+B_{1}B_{3}\right) \nonumber \\
&\quad \quad \quad - \ B_{3}\left(B_{2}B_{1}+B_{3}T_{3}\right), \\
&T_{1} = b_{22}+b_{33}, \ \ T_{2}=b_{11}+b_{33}, \ \ T_{3}=b_{11}+b_{22},\\
&B_{3} = b_{12}, \ \  B_{2}=b_{13}, \ \  B_{1}=b_{23},
\end{align}
and $w_{1},w_{2}$, and $w_{3}$ are given in (\ref{w:1})-(\ref{w:3}).





\bibliographystyle{unsrt}
\bibliography{RootConfV7}



\end{document}